\documentclass[letterpaper, 10 pt, conference]{ieeeconf}  % Comment this line out if you need a4paper
\usepackage{cite, graphicx, amsmath, amssymb, mathrsfs, float, subfig}
\usepackage[font={small}]{caption}
\renewcommand{\baselinestretch}{0.96}

\IEEEoverridecommandlockouts                              % This command is only needed if 
\title{\LARGE \bf
Vehicle Routing Problem with Time Windows: A Deterministic Annealing approach}

\author{Mayank Baranwal$^{1,a}$, Pratik M. Parekh$^{1,b}$, Lavanya Marla$^{2,c}$, Srinivasa M. Salapaka$^{1,d}$, Carolyn L. Beck$^{2,e}$% <-this % stops a space
\thanks{$^{1}$Department of Mechanical Science and Engineering, University of Illinois at Urbana-Champaign, 61801 IL, USA}%
\thanks{$^{2}$Department of Industrial and Enterprise Systems Engineering, University of Illinois at Urbana-Champaign, 61801 IL, USA}%
\thanks{$^{a}${\tt\small baranwa2@illinois.edu}, $^{b}${\tt\small pparekh2@illinois.edu}}
\thanks{$^{c}${\tt\small lavanyam@illinois.edu}, $^{d}${\tt\small salapaka@illinois.edu}}
\thanks{$^{e}${\tt\small beck3@illinois.edu}}
}

\begin{document}

\maketitle
\thispagestyle{empty}
\pagestyle{empty}

%%%%%%%%%%%%%%%%%%%%%%%%%%%%%%%%%%%%%%%%%%%%%%%%%%%%%%%%%%%%%%%%%%%%%%%%%%%%%%%%
\begin{abstract}
The Vehicle Routing Problem with Time-Windows (VRPTW) is an important problem in allocating resources on networks in time and space. We present in this paper a Deterministic Annealing (DA)-based approach to solving the VRPTW with its aspects of routing and scheduling, as well as to model additional constraints of heterogeneous vehicles and shipments. This is the first time, to our knowledge, that a DA approach has been used for problems in the class of the VRPTW. We describe how the DA approach can be adapted to generate an effective heuristic approach to the VRPTW. Our DA approach is also designed to not get trapped in local minima, and demonstrates less sensitivity to initial solutions. The algorithm trades off routing and scheduling in an $n$-dimensional space using a tunable parameter that allows us to generate qualitatively good solutions. These solutions differ in the degree of intersection of the routes, making the case for transfer points where shipments can be exchanged. Simulation results on randomly generated instances show that the constraints are respected and demonstrate near optimal results (when verifiable) in terms of schedules and tour length of individual tours in each solution.
\end{abstract}

%%%%%%%%%%%%%%%%%%%%%%%%%%%%%%%%%%%%%%%%%%%%%%%%%%%%%%%%%%%%%%%%%%%%%%%%%%%%%%%%
\section{INTRODUCTION}

In this paper, our focus is on a class of network-based resource allocation problems that are commonly referred to in the optimization literature as Vehicle Routing Problems with time-Windows (VRPTW). These problems involve both routing and scheduling elements, and are NP-hard in general. VRPTW problems involve the routing and scheduling of multiple vehicles from a depot (or multiple depots) to meet demands at multiple locations, under some time-window considerations. This class of problems is at the core of several problems that involve routing and/or scheduling, and occurs in several application domains including, but not limited to, transportation, logistics and communications. Specific examples include airline scheduling, vehicle routing, service network design, load distribution, production planning, computer scheduling, portfolio selection, and apportionment. These application domains have important economic value, and high importance is attached to achieving economically and computationally efficient solutions to the VRPTW.

VRPTW problems have been some of the most extensively studied problems in the optimization literature. In the past three decades, approaches including exact methods, heuristics, and neighborhood search approaches, or a combination of these have been used. As such, several survey papers on the topic exist. We refer the reader to \cite{VehicleRoutingSurvey}, \cite{TranspOnDemand}, \cite{VRP_book} for the most recent surveys on general vehicle routing problems, and to [1] for some special cases of the VRPTW. Because these problems involve the allocation of multiple heterogeneous and discrete resources, have non-linear costs and constraints, the problems are typically non-linear and non-convex. Heuristic approaches are often subject to achieving local minima rather than global minima, and therefore are dependent on the starting solution that is used. 

Specifically, we consider $K$ vehicles (resources) that start at a depot, and need to pick up shipments from various customer locations, denoted $i = 1,...,N$. Each vehicle $j$ can potentially have capacity $W_j$. Each customer has a time-window [$t_{i,\text{start}},t_{i,\text{end}}$] within which the shipment is to be picked up. Associated with each customer/shipment $i$ is a priority weight $p_i$. Additionally, we are also given the travel times between each pair of points in the network including the depot and customer locations. Our goal is to find the set of customers to be visited by each individual vehicle, the sequence and the schedule according to which they should be visited, while minimizing the costs associated with shipments served earlier or later than their time-windows, and the total distance covered by the vehicles. We also consider the case of heterogeneous shipments, where each shipment can also be one of $l=1,...,p$ types. Each vehicle $j$, has a holding capacity for the $l$th type of shipment, that we denote as $W_{jl}$.

In this paper we present a solution methodology based on Deterministic Annealing (DA) for the VRPTW. DA is well-suited to combinatorial  clustering/resource allocation problems that require obtaining an optimal partition of an underlying domain, and optimally assigning resources to each cell of the partition. DA-based methods have been reported in a vast number of applications such as minimum distortion problems in data compression \cite{rose1998deterministic}, model aggregation \cite{xu2014aggregation}, routing problems in multiagent networks \cite{kale2012maximum}, locational optimization problems \cite{salapaka2003constraints}, and coverage control problems \cite{xu2014clustering}.

Such problems are highly non-convex, computationally complex and suffer from poor local minima that riddle the cost surface \cite{gray1982multiple}. A variety of heuristic approaches have been proposed to address these difficulties, and they range from repeated optimization with different initialization, and heuristics to good initialization, to heuristic rules for cluster splits and merges. In this context, simulated annealing algorithm \cite{kirkpatrick1983optimization}, which capitalizes on the analogy of annealing process in physical chemistry, was shown to be a good iterative method that would achieve a global minimum, but with an annealing rate so slow that this algorithm loses practicality in many applications. In \cite{braekers2011deterministic, braysy2008effective}, a threshold accepting deterministic variant of the simulated annealing algorithm has been proposed. It should be remarked that in some existing literature, such deterministic variants are unfortunately also referred to as the deterministic annealing algorithms; however they are significantly different from the maximum-entropy principle based deterministic algorithm \cite{rose1990deterministic} used in this article - in terms of the underlying heuristics, implementation, goals, and performance. The DA algorithm used in this paper has the ability to avoid poor local optima while still maintaining a relatively fast convergence rate. However, to the best of our knowledge, DA has not been availed to address the specifications and constraints  that arise in the context of simultaneous routing and scheduling problems. Our approach adapts these specifications and constraints in terms of facility location/clustering, thereby enabling the use of DA. Since DA allows for incorporating multiple constraints and different types of specifications, it is well suited to VRPTW (including cases when complex additional constraints are present). For instance, this approach offers the flexibility to tradeoff routing and scheduling in space and time-dimensions to identify clusters in the combined space to allow us to tradeoff between routing and scheduling decisions. We discuss this in detail in the paper.
 
The main contributions of this paper are to incorporate the DA method to VRPTW, where specific capacity, routing, and scheduling constraints are incorporated simultaneously. Simulation results show that the constraints are respected and demonstrate near optimal results (when verifiable). This method provides a tunable parameter, whereby the relative importance of scheduling vs routing is controlled.  Studies with respect to this parameter evidence general features, which can be used for updating scheduling and routing procedures; for instance we show that when schedule is emphasized over routes, there is an increased overlapping of routes and self  intersections, which suggests that such scenarios will benefit by having {\em transfer/exchange locations}  where shipments between different vehicle routes can be  exchanged.

\section{PROBLEM DESCRIPTION}\label{sec:problem_description}
In order to solve the VRPTW, we first describe a series of related problems, which are in increasing order of complexity.

\subsection{Uncapacitated Scheduling Problem}\label{subsec:unCap-SP}
This problem involves finding a schedule according to which each vehicle will visit the customer (shipment) locations in order to pick up shipments. In this problem we ignore the locations of customers and assume travel between any two locations is instantaneous. While this is not applicable in practice, this helps illustrate the working of our algorithm. Note that in the context of VRPTW, this problem only considers the vehicle allocation for a given schedule information; it does not have any routing aspect to it. More precisely, the problem reduces to a scenario where shipments with prescribed time-windows are needed to be picked up from a depot; the objective is to  ascribe an arrival time for each vehicle so that the maximum number of shipments are serviced. 

\subsection{Capacitated Scheduling Problem}\label{subsec:Cap-SP}
This problem builds upon the uncapacitated scheduling problem by adding capacity constraints on the number of customers (shipments) that can be carried by each vehicle. Vehicles are heterogeneous, that is, can have varying capacities. 

\subsection{Scheduling Problem with Multiple Capacity Constraints}\label{subsec:CapMulti-SP}
In this version of the capacitated scheduling problem, customers (shipments) picked up have varying types. They thus have different occupancy rates in each vehicle. Each vehicle has a capacity associated with each shipment type it can carry.

\subsection{Vehicle Routing Problem (VRP)}\label{subsec:VRP}
In the VRP, we find a set of locations that is served by each vehicle in the fleet, as well as the sequence according to which each customer location should be visited. We assume that there are no time-window restrictions on when each customer (shipment) can be visited/served.

\subsection{Vehicle Routing Problem with Time Windows (VRPTW)}\label{subsec:VRPTW}
The VRPTW involves solving the VRP with additional time-window restrictions on the time that each customer can be visited. This problem can incorporate heterogeneous shipment sizes as well as heterogeneous capacities. 

\section{DETERMINISTIC ANNEALING ALGORITHM: A MAXIMUM ENTROPY PRINCIPLE APPROACH FOR CLUSTERING}\label{sec:DA_intro}
At its core, the DA algorithm solves a facility location problem (FLP): For given $N$ customer locations, find $K$ facility locations such that the total {\em weighted sum of the distance of each customer to its nearest facility is minimized}. In other words, if $s_i$ and $r_j\in\mathbb{R}^n$ denote the locations of $i^{th}$ customer and $j^{th}$ facility, respectively, then the FLP addresses the following optimization problem: 
\begin{equation}\label{eq:problem_def}
	\min\limits_{r_j\in\Omega,1\leq j\leq K}\sum\limits_{i=1}^{N}p_i\left\{\min\limits_{r_j,1\leq j\leq K}d(s_i,r_j)\right\}
\end{equation}
where, $d(s_i,r_j)\in\mathbb{R}_+$ denotes the distance between the $i^{th}$ customer location $s_i$ and $j^{th}$ facility location $r_j$, $\Omega\in\mathbb{R}^n$ is a compact domain of interest and $p_i$ is a given positive constant (without loss of generality, we assume $\sum_ip_i = 1$) that denotes the relative weight of the $i^{th}$ customer. Note that in general, $s_i$ and $r_j$ need not necessarily denote physical locations but can belong to any relevant {\em property space} $\Omega$. Borrowing from the data compression literature \cite{cover2012elements}, we define {\em distortion} as a measure of the average distance of a customer to its nearest facility, given by $D(s,r) = \sum\limits_{i=1}^{N}p_i\min\limits_{1\leq j\leq K}d(s_i,r_j)$. Then the equivalent optimization problem is to minimize the distortion function. The solution to a FLP essentially results in a series of clusters, where the facility $j$ is located at the centroid $r_j$ of the $j^{th}$ cluster and each customer is associated only to its nearest facility.

Note that any change in location of a particular customer $i$ affects $d(s_i,r_j)$ only with respect to the {\em nearest} facility $j$ in the distortion. This distributed aspect makes most algorithms (such as Lloyd's \cite{lloyd1982least}) overly sensitive to the initial facility location. The DA algorithm \cite{rose1990deterministic} overcomes this sensitivity by modifying the distortion $D$ such that every customer $i$ is associated to every facility $j$ through an {\em association weight} $p(j|i)$:
\begin{small}
\begin{equation}\label{eq:distortion}
	\bar{D}(s,r) = \sum\limits_{i=1}^{N}p_i\sum\limits_{j=1}^{K}p(j|i)d(s_i,r_j).
\end{equation}
\end{small}
Here we choose $\{p(j|i)\}$ to satisfy $0\leq p(j|i)\leq 1$ and $\sum\limits_{j=1}^{K}p(j|i)=1$ without loss of generality. Thus, we replace the {\em average distance of a customer to its nearest facility} by the {\em weighted average distance of a customer to all the facilities}. The weights $p(j|i)$ assess the trade-off between decreasing the {\em local} influence and the deviation of $\bar{D}$ from the original distortion $D$. Since we have no prior knowledge on the association weights, we apply the principle of maximum entropy to estimate them. The Shannon entropy term $H(r|s) = -\sum\limits_{i=1}^{N}p_i\sum\limits_{j=1}^{K}p(j|i)\log(p(j|i))$, widely used in data compression literature \cite{cover2012elements}, measures uncertainties in facility locations with respect to the known customer locations. Thus, maximizing the entropy is equivalent to decreasing the {\em local} influence.

This trade-off between minimizing the distortion in (\ref{eq:distortion}) and maximizing the entropy \cite{rose1998deterministic, rose1990deterministic} is achieved by seeking $\{p(j|i)\}$ that minimize the {\em free energy}, or the Lagrangian, given by $F(r)\triangleq \bar{D}(s,r)-TH(r|s)$, where $T$ is a Lagrange multiplier, referred to as {\em temperature}. This yields a {\em Gibbs distribution}
\begin{small}
\begin{equation}\label{eq:gibbs}
	p(j|i) = \frac{\exp\{-\beta d(s_i,r_j)\}}{\sum\limits_{k=1}^{K}\exp\{-\beta d(s_i,r_k)\}}
\end{equation}
\end{small}
with $\beta = 1/T$. By substituting the association weights (\ref{eq:gibbs}), the free energy in (\ref{eq:free_energy}) simplifies as
\begin{small}
\begin{equation}\label{eq:free_energy}
	F(r) = -\frac{1}{\beta}\sum\limits_{i=1}^{N}p_i\log\sum\limits_{k=1}^{K}\exp\{-\beta d(s_i,r_k)\}.
\end{equation}
\end{small}
In the DA algorithm,  this {\em free energy} function is then deterministically optimized at successively reduced temperatures over repeated iterations.

The readers are encouraged to refer to \cite{parekh2015deterministic} for detailed analysis on the complexity of the DA algorithm. For implementation on very large datasets, a scalable modification of the DA is proposed in \cite{sharma2006scalable}.

\section{SOLUTION APPROACH: MODIFICATIONS OF THE DA ALGORITHM}\label{sec:sol_approach}
DA addresses both routing and scheduling problems as modifications of the clustering problem, wherein by choosing $\Omega$ as a domain of spatial coordinates, the DA solution groups together customers served by the same vehicle, while clustering on the space $\Omega$ of time-windows results in service schedules for each customer. We now describe the DA algorithm in the context of VRP to accommodate various scheduling, routing and capacity constraints.

\subsection{DA for scheduling}\label{subsec:DA_schedule}
In this scenario, the DA is used to allocate vehicle arrival times $r_j$ to service the shipments with specified service time-windows $[t_{i,\text{start}}, t_{i,\text{end}}]$. Accordingly we choose $\Omega\subset\mathbb{R}$ representing the time domain of interest, the property $s_i$ of the customer $i$ is chosen to be the {\em mid-point} of the associated time window, i.e.,
\begin{small}
\begin{equation}\label{eq:si_schedule}
	s_i = \frac{t_{i,\text{start}} + t_{i,\text{end}}}{2};
\end{equation}
\end{small}
a convenient choice of the distance metric between the customer $i$ and facility $j$ is the squared-euclidean distance $d(s_i,r_j) = |s_i-r_j|^2$. The distortion $\bar{D}$ captures the penalty of deviating from the {\em mid-times}. Thus, minimizing the free-energy $F$ in (\ref{eq:free_energy}) is commensurate with the cost incurred for not serving within the specified time-window being minimized over all shipments.

\subsection{DA for routing}\label{subsec:DA_routing}
Some common heuristics to the standard single depot vehicle routing problems include approaches, such as, {\em cluster first-route second} \cite{gillett1974heuristic}. We follow the same approach wherein we first cluster the demand nodes based on their geographical distances, and then later design economical routes over each cluster. Executing the DA algorithm directly on the space of geographical coordinates results in a partition of the underlying domain $\Omega$, wherein each cluster $j$ is served by a different vehicle. In the clustering process $\Omega\subset\mathbb{R}^2$ is the spatial domain of interest and we use the parameter $\beta$ to cluster locations. The heuristic behind this approach is that if we use $K=N$ resources and let $\beta\to\infty$ then we expect each customer to be associated with one resource (vehicle).

However, the clustering lacks the sequencing aspect of the routing problem. In order to incorporate the sequencing aspect, we include the minimum tour-length constraint in the original formulation, which essentially amounts to solving a traveling salesman problem (TSP) in each cluster. This is achieved by appropriately modifying the free energy in (\ref{eq:free_energy}) to obtain the new Lagrangian
\begin{small}
\begin{equation}\label{eq:lagrangian_routing}
	F'(r) = F(r) + \lambda\left(\sum\limits_{k=1}^{K}d(r_j,r_{j-1})-L\right),
\end{equation}
\end{small}
where $L$ is a given tour-length, and $\lambda$ is the Lagrange multiplier related to it. As $L$ is gradually increased from a small value, the tour length between locations $\{r_j\}$  in the same cluster would coincide with  the optimal tour length between $\{s_i\}$. Here the spatial coordinates $\underbar{x}_i = (x_i,y_i)$ represent the property $s_i$ of the customer $i$. The distance function $d(s_i,r_j)$ between the customer $i$ and facility $j$ is chosen as : $d(s_i,r_j) = \|s_i-r_j\|_2^2$. We adopt the Elastic Net approach \cite{durbin1989analysis} for finding the optimal tour-length. The required tour-length is controlled by appropriately varying the Lagrange multiplier $\lambda$ with $\beta$ \cite{rose1990deterministic, roehl2011maximum}.

\subsection{DA for scheduling-cum-routing}\label{subsec:DA_scheduling_cum_routing}
In this setting, we desire a resource $j$ to be located at a customer $i$ in a pre-ascribed time-window. In this context, the property $s_i$ of the customer $i$ is given by,
\begin{small}
\begin{equation}\label{eq:si_schedule_cum_route}
	s_i = \left( \begin{array}{c}
	x_i\\
	y_i\\
	\frac{\lambda}{2}(t_{i,\text{start}}+t_{i,\text{end}}) \end{array} \right),
\end{equation}
\end{small}
where $\underbar{x}_i = (x_i,y_i)$ and $[t_{i,\text{start}}\quad t_{i,\text{end}}]$ represent the spatial coordinates and service time-window of the customer $i$, respectively, and the domain of interest is $\Omega\subset\mathbb{R}^2\times\mathbb{R}_+$. We define a new metric that incorpates both {\em scheduling} and {\em routing} through a {\em velocity} parameter $\lambda$ 
\begin{small}
\begin{equation}\label{eq:new_metric}
	d(s_i,r_j) = \left\|\underbar{x}_i-\begin{bmatrix} r_j(1)\cr r_j(2)\end{bmatrix}\right\|_2^2 + \frac{\lambda^2}{4}|(t_{i,\text{start}}+t_{i,\text{end}})-r_j(3)|^2.
\end{equation}
\end{small}
Clearly, setting $\lambda = 0$ solves the routing problem, as described in Sec. \ref{subsec:DA_routing}, while setting $\lambda$ to a very high value virtually solves the scheduling problem in Sec. \ref{subsec:DA_schedule}. Thus $\lambda$ captures the trade-off between {\em distance} based optimization and {\em time-window} based optimization problems. In practice, $\lambda$ should be set to a value comparable to the average speed of the vehicles.

\subsection{DA for heterogeneous capacitated resources/vehicles}\label{subsec:single_capacity}
In this setting, we consider vehicles $\{j\}$ to be of different types. If $\lambda_j$ denotes the total allocated capacity of the resource $j$, i.e. if $\lambda_j$ amount of resources are located at $r_j$, then the modified Gibbs distribution is given by $p(j|i) = \frac{\lambda_j\exp\{-\beta d(s_i,r_j)\}}{\sum\limits_{k=1}^{K}\lambda_k\exp\{-\beta d(s_i,r_k)\}}$. Note that the parameter $\lambda_j$ in the modified distribution specifies the weight of the $j$th resource; the resulting solution $p_j:p_k$ tend to be in the same ratio as $\lambda_j:\lambda_k$, where $p_j$ is marginal distribution of the $j$ resource given by $p_j=\sum_i p_ip(j|i)$ \cite{rose1998deterministic, salapaka2003constraints, salapaka2003locational}.

\subsection{DA for heterogeneous capacitated vehicles and heterogeneous shipments}\label{subsec:multiple_capacity}
In this setting, we have $p$ types of customers, and the total amount of resource $j$ needing to be allocated to all customers of $l^{th}$-type are $\lambda_{jl}, 1\leq l\leq p, 1\leq k\leq K$. These constraints are incorporated in the Gibbs distribution as $p(j|i) = \frac{\sum\limits_{l=1}^{p}\lambda_{j,l}\exp\{-\beta d(s_i,r_j)\}}{\sum\limits_{k=1}^{K}\sum\limits_{l=1}^{p}\lambda_{kl}\exp\{-\beta d(s_i,r_k)\}}$.
\begin{figure*}
	\begin{center}
	\begin{tabular}{ccc}
	\includegraphics[width=0.79\columnwidth]{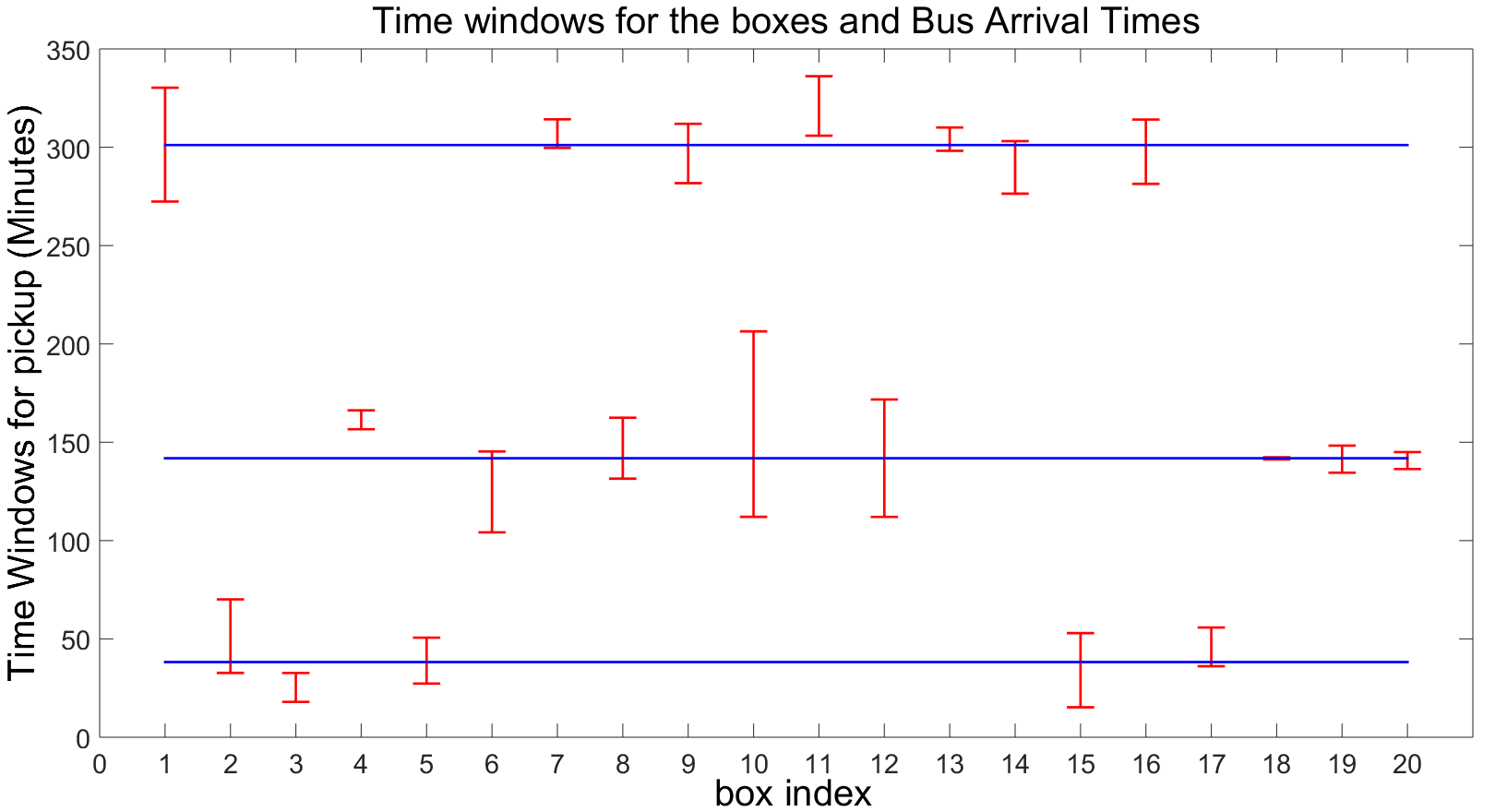}&\includegraphics[width=0.79\columnwidth]{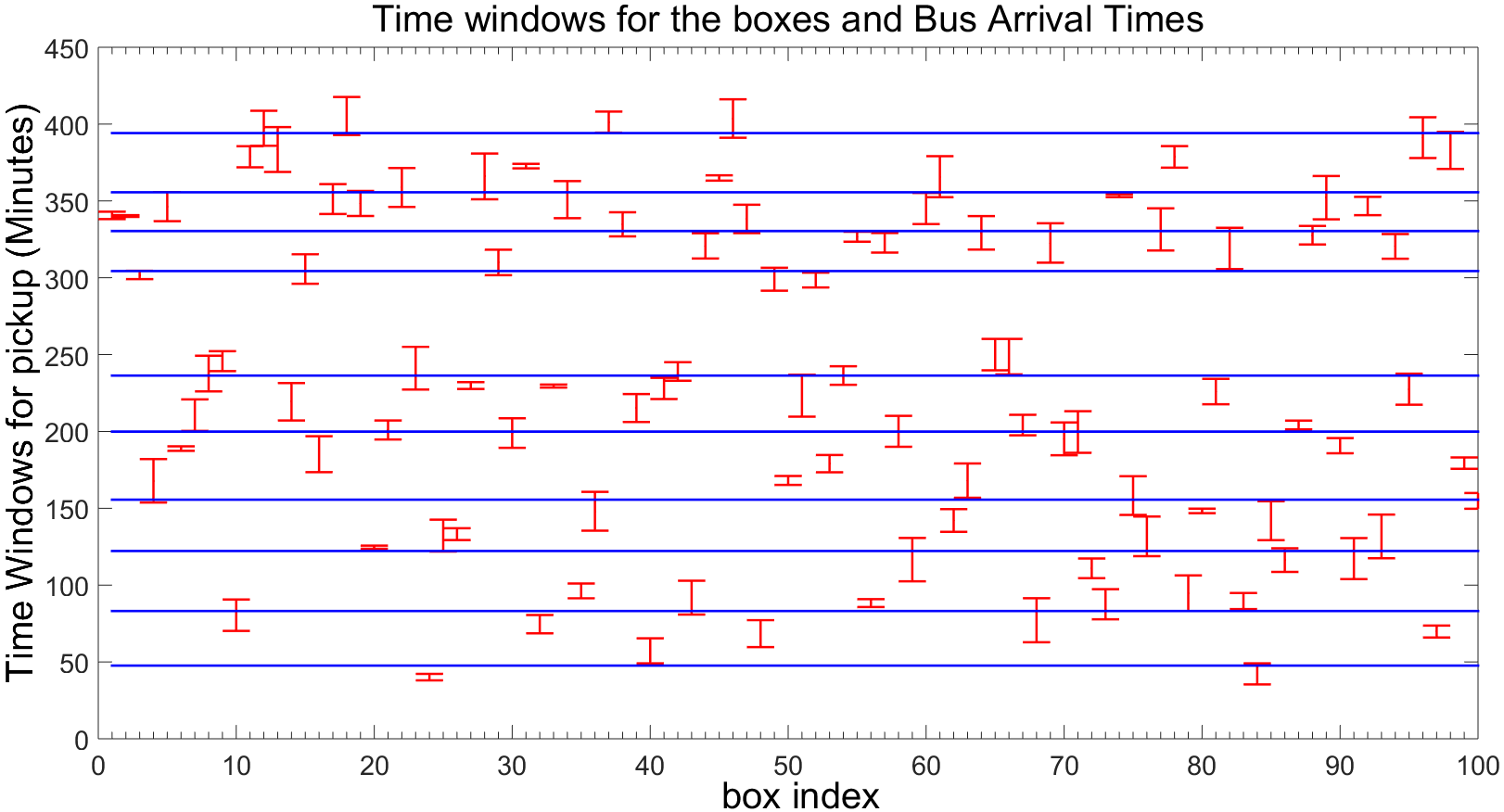}&\includegraphics[width=0.34\columnwidth]{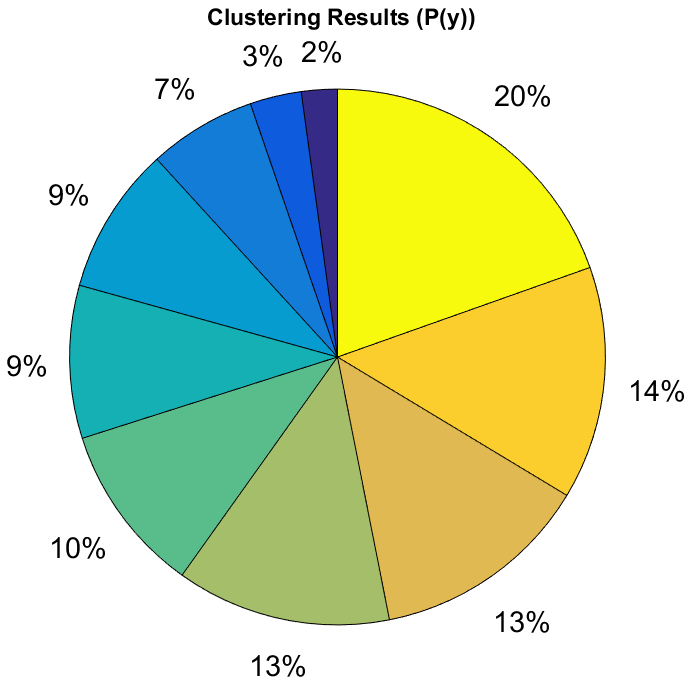}\cr
	(a)&(b)&(c)\end{tabular}
	\vspace{-0.5em}
	\caption{(a) Time-Windows for 20 shipments (Red Intervals) and the Arrival Times (Cluster locations (blue lines)) for 3 vehicles. (b) Time-Windows for 100 shipments (Red Intervals) and the Arrival Times (Cluster locations (blue lines)) for 10 vehicles. (c) Capacity associated with each cluster for the 100 shipments scenario. Clearly, the capacities obtained by the DA algorithm are in proportion to the capacity values given as constraints.}
	\label{fig:P1}
	\vspace{-1em}
	\end{center}
\end{figure*}

\section{SIMULATIONS AND RESULTS}\label{sec:sol}
In this section, we demonstrate the effectiveness of DA for solving scheduling and routing problems of varying degrees of complexity. We generate randomized instances of each problem type and the algorithm is tested for each of the problem instances. However, for the sake of demonstration, we show only one instance of each problem type.

\subsection{Uncapacitated Scheduling Problem}\label{subsec:result_unCap-SP}
For the uncapacitated scheduling problem, we consider an instance with $20$ shipments needing to be serviced by $3$ vehicles. The depot incurs the corresponding cost if the shipment does not get serviced. That is, if a vehicle arrival time does not lie in the time-window for a shipment, the shipment goes to waste and the corresponding cost is incurred by the depot. We randomly sample the time-windows to lie within $[0,350]$ minutes.

In Fig. \ref{fig:P1}a, blue lines denote the final clustered information (vehicle arrival times) after direct execution of the DA algorithm (as described in Sec. \ref{subsec:DA_schedule}). The red bars denote the time-windows for each shipment. Only three shipments indexed $3, 4$ and $11$ do not have any vehicle arrivals within their time-windows. From the time-intervals seen in Fig. (\ref{fig:P1}a), it is evident that the omission of these shipments is inevitable, as accommodating any of them leads to larger costs as other shipments will not be served.
%{\bf Remark:} The DA algorithm was also compared with the Lloyd's algorithm. Simulations were performed on datasets with $20$ shipments and the vehicle arrival times $(3)$ were computed using Lloyd's algorithm. The DA algorithm took relatively greater computational time on average than the Lloyd's algorithm, however the Lloyd's algorithm was more prone to getting trapped in local minima; moreover, the solutions were highly sensitive to initialization. It should be noted that the number of iterations over $\beta$ is not large, since the parameter $\beta$ is increased geometrically from $\beta$ = $\epsilon$ to $\beta_\text{stop}$. Also, each iteration in the DA algorithm, as $\beta$ tends to infinity can be shown to have the same computational complexity as the Lloyd's algorithm. Therefore it is expected that the computational complexity for the scalable DA algorithm is roughly equivalent to a finite (and small) number of iterations of the Lloyd's algorithm.

\subsection{Capacitated Scheduling Problem}\label{subsec:result_Cap-SP}
In the capacitated scheduling problem, we provide capacity to each vehicle and the algorithm should attempt to match these capacities. Once again, the time-windows are randomly sampled to lie within $[0,350]$ minutes. We have $100$ shipments needing to be picked up by only 10 vehicles, where the capacities assigned to each vehicle (cluster) are $\{3,5,9,9,9,10,13,13,13,16\}$. 

Subjected to these capacities, the algorithm works well and only $30\%$ of the shipments are not serviced. The effectiveness of the DA algorithm for clustering with capacity constraints is observed in the pie chart in Fig. \ref{fig:P1}c. The pie chart shows the fraction of shipments associated with each vehicle, which is mathematically depicted by $p(r_j) = \sum_ip_ip(j|i)$. The capacities associated with each cluster (vehicle) are in proportion with the capacity values given as constraints.

\subsection{Scheduling Problem with Multiple Capacity Constraints}\label{subsec:result_CapMulti-SP}
For this problem, there are capacities associated with each vehicle and all types of shipments. For the simulation, we choose an instance with $100$ shipments in total, $3$ types of shipments and only $10$ vehicles for pickup. There are $34$ shipments of type-$1$ (red), $36$ shipments of type-$2$ (magenta), and $30$ shipments of type-$3$ (green). The capacities are randomly chosen and subject to these capacities, the algorithm works well. 

\begin{figure*}
	\begin{center}
	\begin{tabular}{ccc}
	\includegraphics[width=0.65\columnwidth]{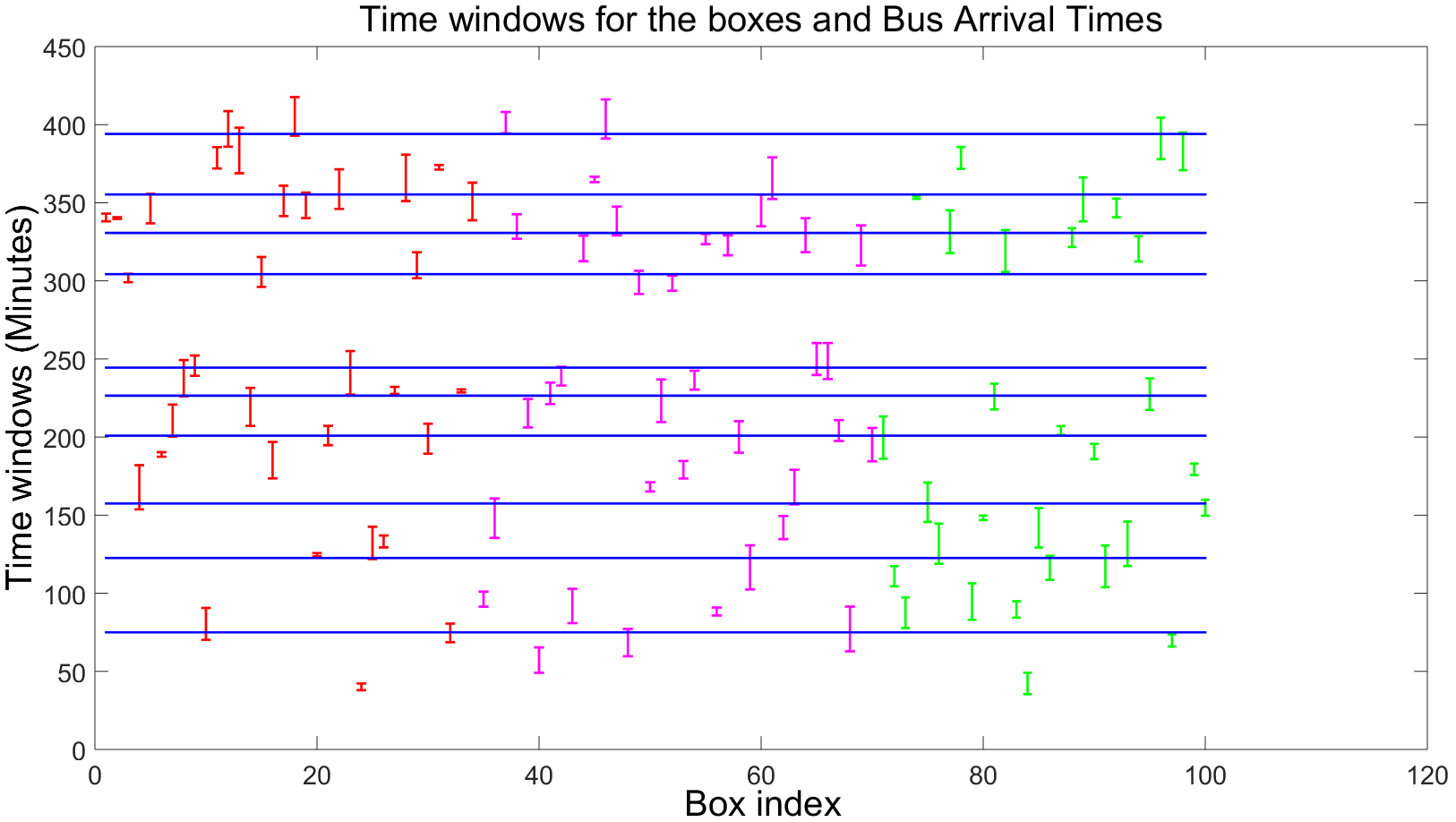}&\includegraphics[width=0.6\columnwidth]{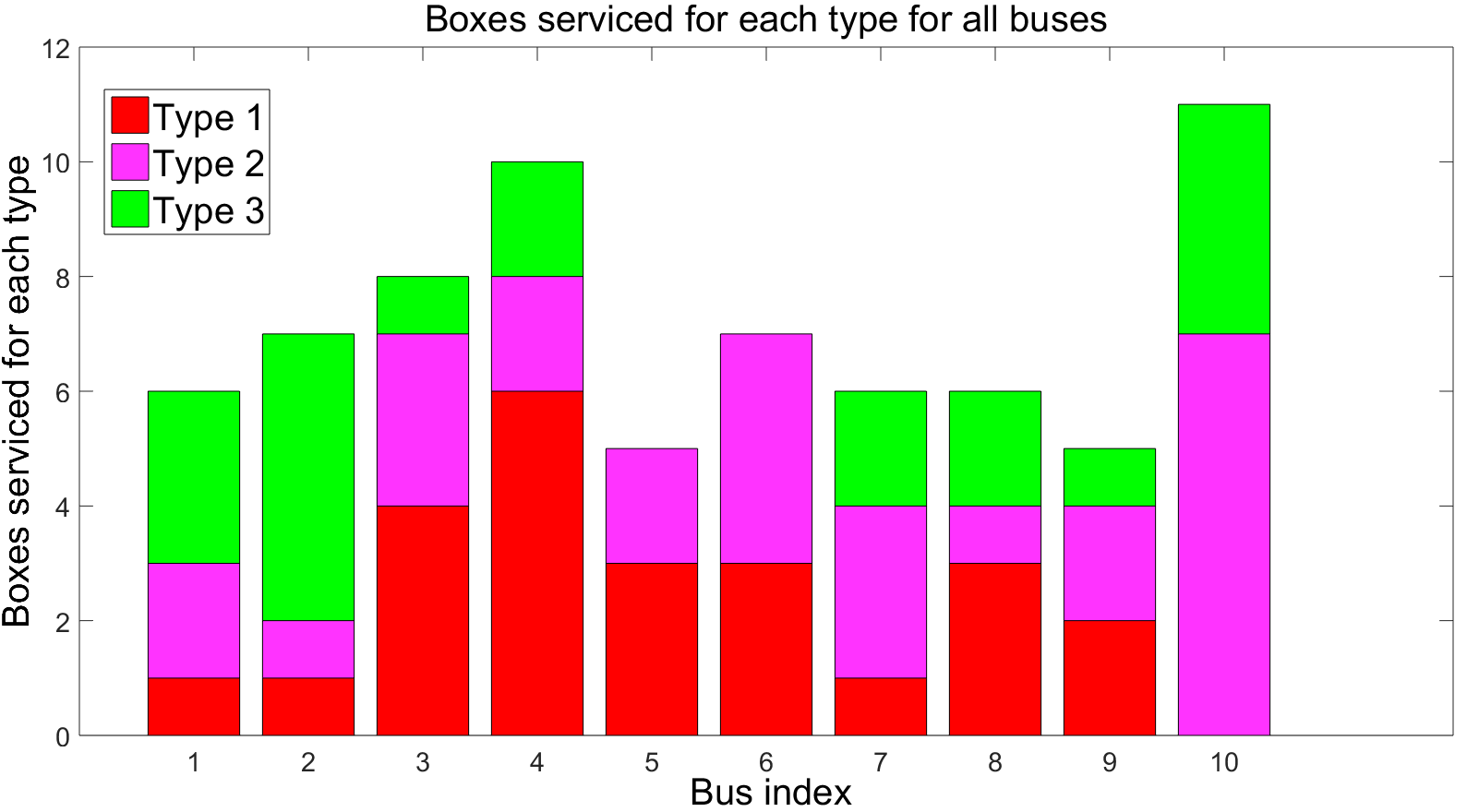}&\includegraphics[width=0.7\columnwidth]{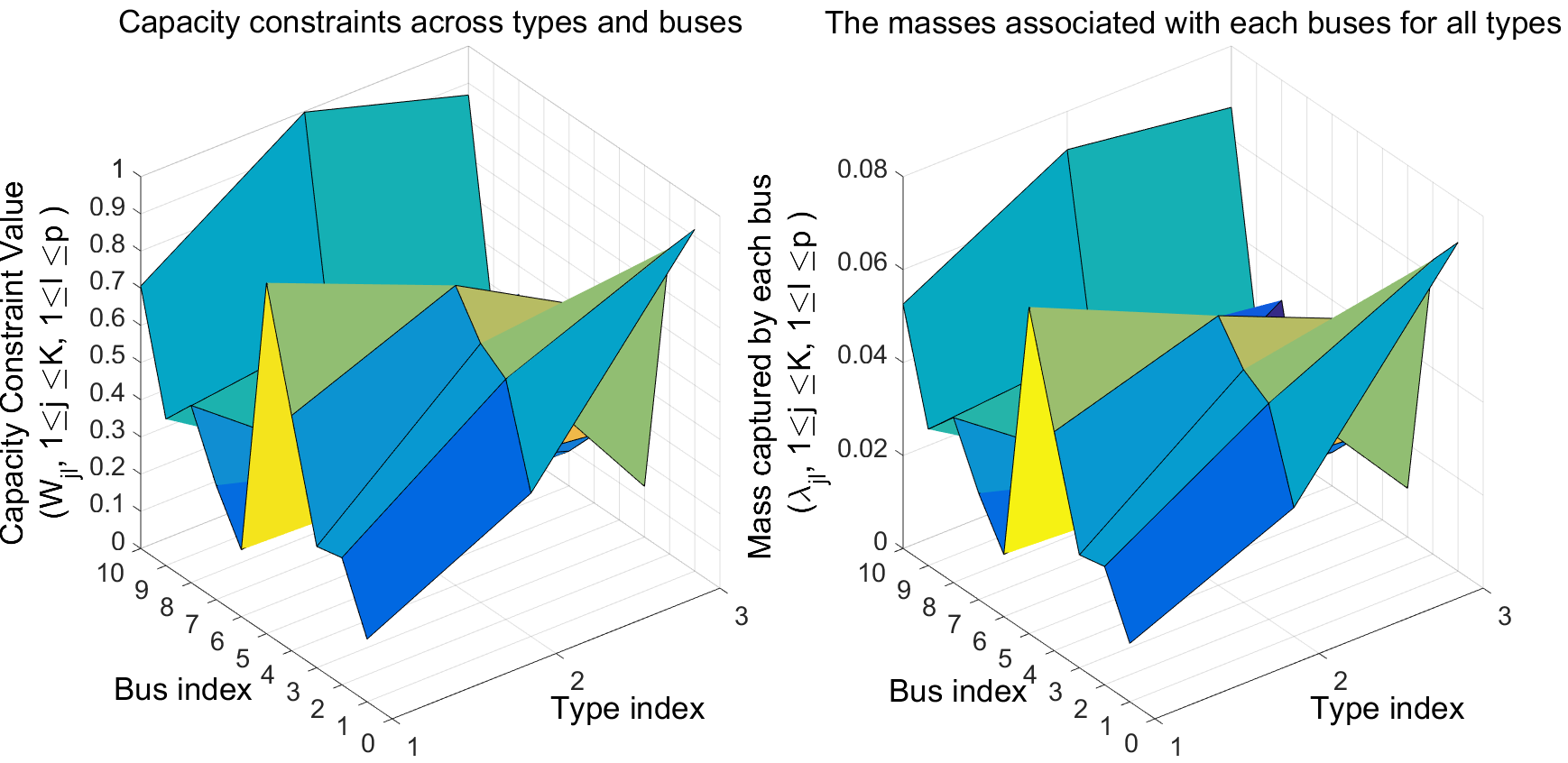}\cr
	(a)&(b)&(c)\end{tabular}
	\vspace{-0.5em}
	\caption{(a) Time-Windows for 100 shipments (3 types) and the Arrival Times (Cluster locations (blue lines)) for 10 vehicles. The three colors (red, magenta and green) are used to indicate the type of shipments. (b) Number of shipments served by each vehicle for each type. (c) [Right]: Capacity constraint values for each type of shipment and vehicle. [Left]: Capacity associated with each vehicle and each type of shipment. Clearly, the shapes indicate that the capacities obtained by the constrained DA algorithm are in proportion with the capacities provided as constraints.}
	\label{fig:P3}
	\vspace{-2em}
	\end{center}
\end{figure*}

In Fig. \ref{fig:P3}a, the blue lines depict the vehicle arrival times, and the red, magenta and green bars indicate the time-windows for the three types of shipments respectively. A total of $71$ shipments are picked up. Fig. \ref{fig:P3}b shows the number of shipments picked up by each vehicle.

Fig. \ref{fig:P3}c shows the effectiveness of the DA algorithm in respecting the capacity constraints. The plot on the left shows the capacity constraints $(\lambda_{jl})$ for each type of shipment and vehicle. The plot on the right shows the clustered mass information. Clearly, the shapes of the two plots match each other.

\subsection{Vehicle Routing Problem}\label{subsec:result-VRP}
For the VRP, we consider a $59$ cities network (and a depot city from where the vehicles start) whose geographical coordinates are located in a $[-100,100]\times[-100,100]$ {unit}$^2$ area. The coordinates of the vehicle depot are given by $(0,0)$. The nodes need to be serviced by a total of six vehicles.

\begin{figure}[tphb]
	\centering
	\includegraphics[width=0.9\columnwidth]{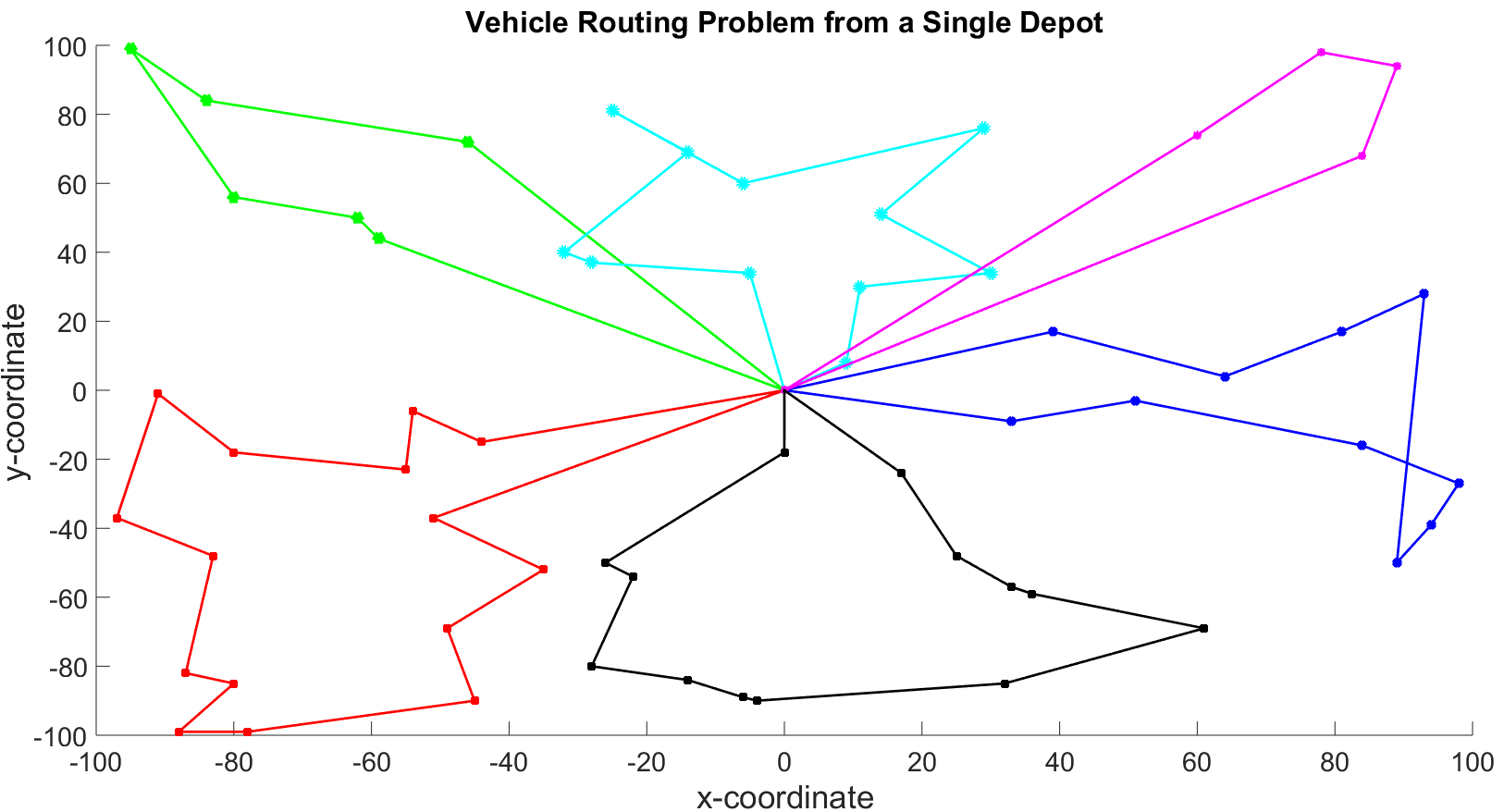}
	\caption{Solution to the VRP using \underline{cluster first-route second} approach. The nodes are first clustered based on their pairwise distances, followed by computation of optimal routes over these clusters.}
	\label{fig:P4}
	\vspace{-1em}
\end{figure}

We first cluster the network based on their pairwise square-Euclidean distance $d(\underbar{x}_i,\underbar{x}_j) = \|\underbar{x}_i-\underbar{x}_j\|_2^2$, i.e., the DA algorithm described in Sec. \ref{subsec:DA_routing} is executed on the set of cities $\{\underbar{x}_i\}$ to obtain $K = 6$ partitions of the entire network. This is followed by solving a TSP in each of the individual clusters using the modification described in Sec. \ref{subsec:DA_routing}. Fig. \ref{fig:P4} shows the effectiveness of the DA algorithm for solving an instance of the VRP.

{\bf Remark}: While the cluster first-route second approach does not necessarily result in an overall global optimal solution, each individual tour obtained by the DA algorithm in the prescribed is indeed optimal (as verified using exhaustive simulations).

\subsection{Vehicle Routing Problem with Time Windows}\label{subsec:result-VRPTW}
For the VRPTW, we again consider the same network of $59$ cities, however, with respective delivery time-windows. The time-windows are chosen to lie in the interval $[0,350]$ minutes. 
The maximum speed of each vehicle is chosen to be $10$ units/min.

\begin{figure*}[!t]
	\begin{center}
	\begin{tabular}{cc}
	\includegraphics[width=0.85\columnwidth]{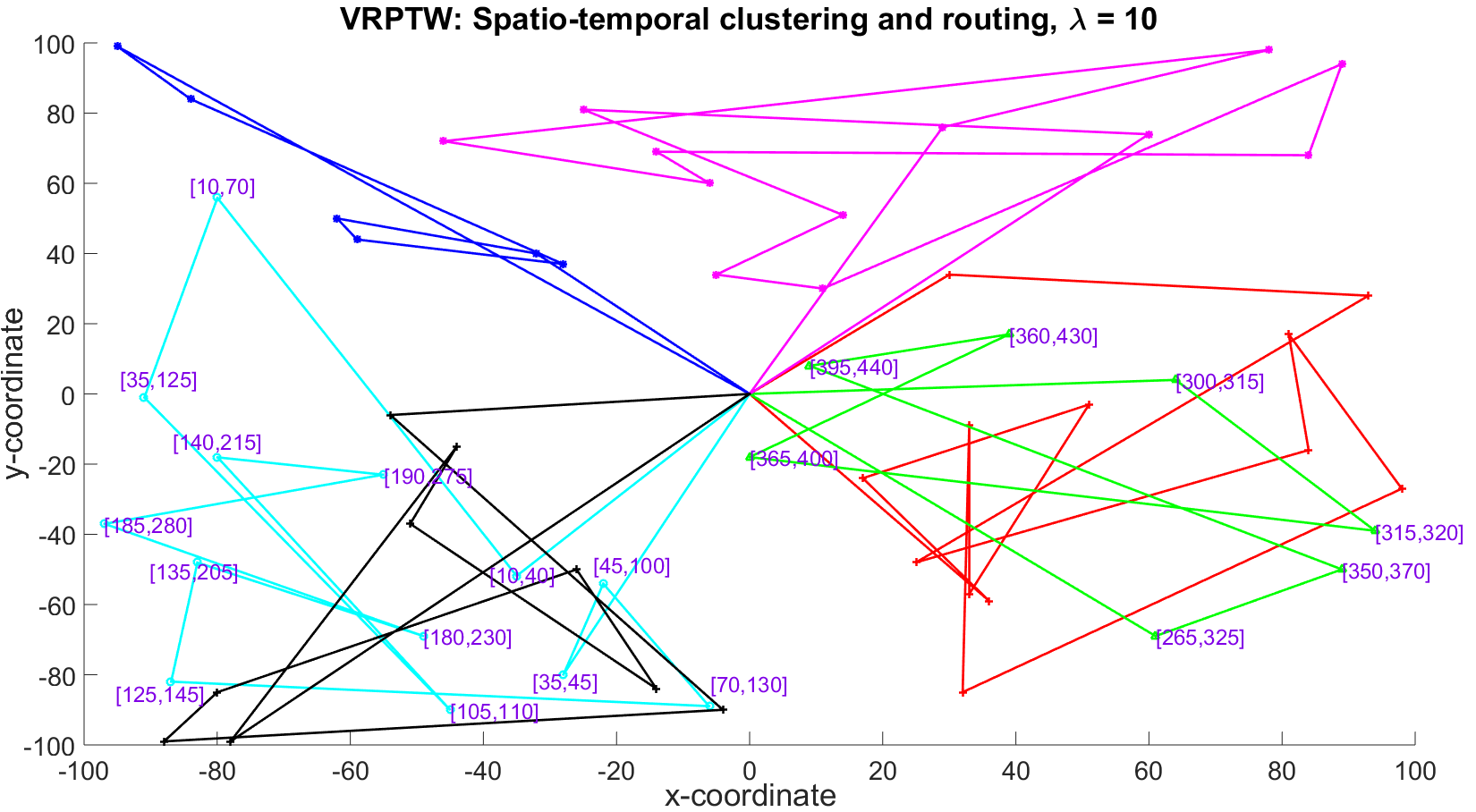}&\includegraphics[width=0.85\columnwidth]{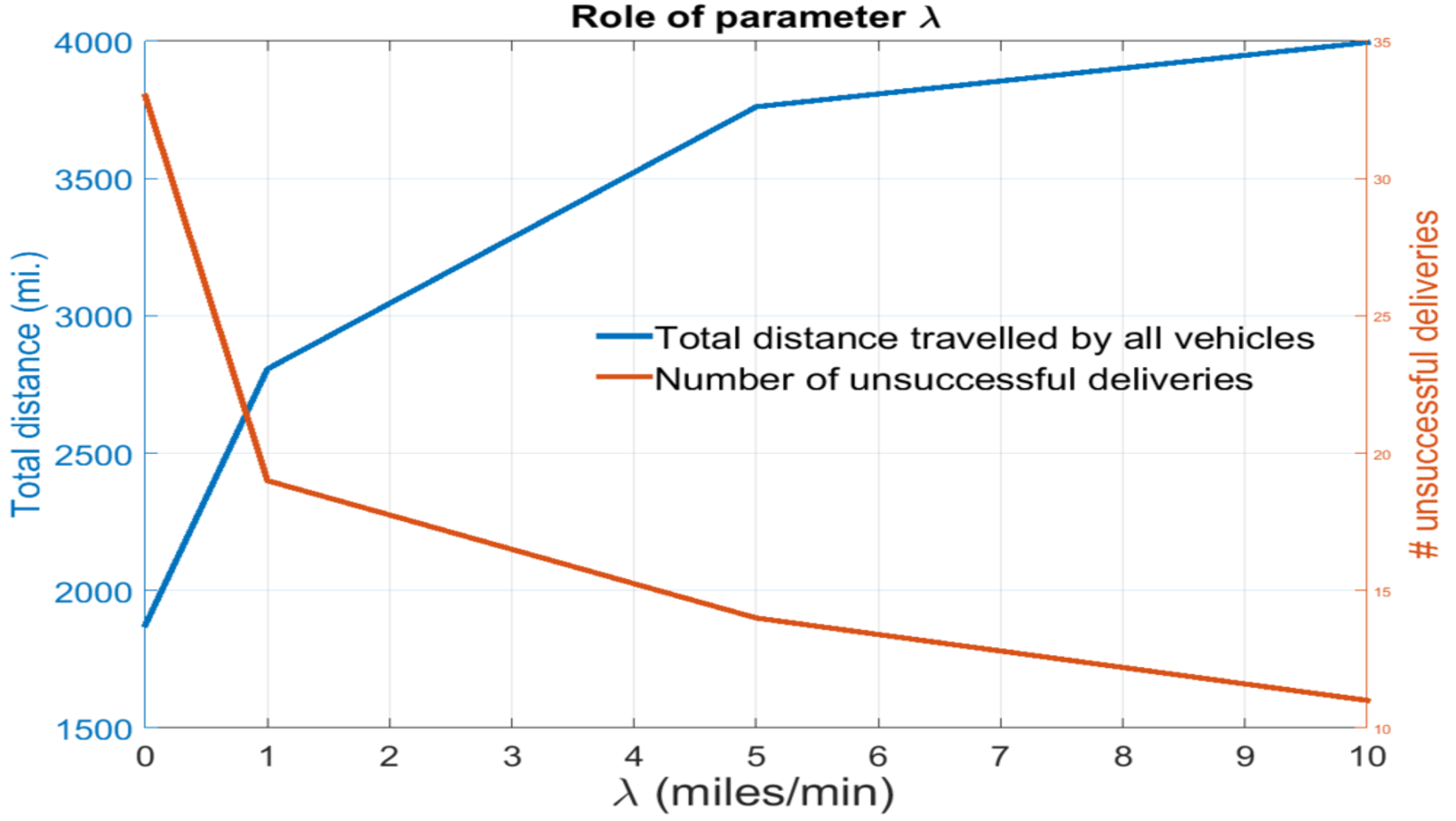}\cr
	(a)&(b)\end{tabular}
	\vspace{-0.5em}
	\caption{(a) Computation of optimal routes by the DA algorithm using the proposed spatio-temporal metric. The method ensures timely delivery of packages so that losses incurred to the depot are minimized. Two such routes (green and cyan) are shown with corresponding time-windows, where the corresponding deliveries for most cities are met on time. (b) Effect of parameter $\lambda$ on the total tour-length and the number of unsuccessful deliveries.}
	\vspace{-2em}
	\label{fig:P5c}
	\end{center}
\end{figure*}

The cities are first clustered using the spatio-temporal metric described in Eq. \ref{eq:new_metric}. Introduction of the additional time coordinate in the problem domain results in intersecting partitions in the spatial coordinates (see Fig. \ref{fig:P5c}), i.e., the paths traversed by different vehicles have overlaps. This observation is in consonance with the fact that vehicles may traverse overlapping routes to ensure timely delivery of packages. Once the cities are clustered in the spatio-temporal domain, the DA algorithm is employed again for solving TSPs over each cluster in this new domain, i.e., the TSP is solved using the new metric. This heuristic ensures timely delivery of shipments so that losses incurred to the depot are minimized. The algorithm achieves this objective by finding self-intersecting (zig-zag) routes. In Fig. \ref{fig:P5c}a, two such routes (green and cyan) are shown with corresponding time-windows. A total of only $11$ shipments are not delivered on time.

It should be noted that the parameter $\lambda$ in (\ref{eq:new_metric}) controls the trade-off between the {\em distance} based optimization and the {\em time-window} based optimization. In other words, $\lambda$ can be interpreted as the average velocity of a vehicle. Smaller values of $\lambda$ result in reduced total tour-lengths, but allow for larger number of not-on-schedule deliveries, as shown in Fig. \ref{fig:P5c}b. As $\lambda$ is increased, the paths traversed by each vehicle overlap and appear to have self-intersections as well as intersections with routes of other vehicles. Thus, such a framework also allows for existence of potential {\em transfer points}, where shipments can be transferred from one vehicle to another. Note that $\lambda = 0$ corresponds to the scenario of VRP (see Fig. \ref{fig:P4}), where routes do not overlap or self-intersect.

\section{CONCLUSIONS AND FUTURE WORK}\label{sec:conclusion}
In this paper, scheduling and routing problems of varying degrees of complexity are addressed in the Deterministic Annealing framework. The algorithm surfaces as an efficient heuristic for solving combinatorially complex scheduling and routing problems. Following are some key features and observations of the DA algorithm (and its modifications):\\
(1) The solutions provided by the DA are completely independent of initializations. In fact, the algorithm attempts to solve a completely deterministic optimization problem (maximization of entropy term) in the beginning.\\
(2) The algorithm has ability to avoid poor local optima. This is illustrated by the quality of solutions for various problem types. For example, in the m-TSP scenario in the VRP in Sec. \ref{subsec:result-VRP}, the individual tour-lengths are found to be optimal.\\
(3) When modified appropriately, the algorithm respects various capacity constraints. For example, in Sec. \ref{subsec:result_CapMulti-SP} we observe that the capacities obtained by the DA algorithm match the required capacity constraints for each shipment and each vehicle.\\
(4) A new spatio-temporal metric is proposed which incorporates the concept of time-based distance using a {\em velocity} parameter $\lambda$. As $\lambda$ is increased, the routes traversed by different vehicles are observed to overlap more, and a relatively larger number of self-intersections are observed within each route to ensure timely delivery of packages. This may facilitate introduction of the {\em transfer points} where shipments from one vehicle are transferred to another to be delivered by the second vehicle.\\
The solutions provided by DA are observed to be robust to uncertainties in customer locations or service schedules, which will be demonstrated in the subsequent work. Another future direction is to design distance functions that incorporate the time-windows directly in the problem formulation. We intend to test our algorithm against the existing heuristics on benchmark instances and the results will be reported in the subsequent work. We also aim to include the notion of {\em transfer points} in our problem formulation as part of our future work.


%%%%%%%%%%%%%%%%%%%%%%%%%%%%%%%%%%%%%%%%%%%%%%%%%%%%%%%%%%%%%%%%%%%%%%%%%%%%%%%%

%%%%%%%%%%%%%%%%%%%%%%%%%%%%%%%%%%%%%%%%%%%%%%%%%%%%%%%%%%%%%%%%%%%%%%%%%%%%%%%%

%%%%%%%%%%%%%%%%%%%%%%%%%%%%%%%%%%%%%%%%%%%%%%%%%%%%%%%%%%%%%%%%%%%%%%%%%%%%%%%%

\section*{ACKNOWLEDGMENT}

The authors would like to acknowledge NSF grants ECCS 15-09302 and CNS 15-44635 for supporting this work.

%%%%%%%%%%%%%%%%%%%%%%%%%%%%%%%%%%%%%%%%%%%%%%%%%%%%%%%%%%%%%%%%%%%%%%%%%%%%%%%%

\bibliographystyle{IEEEtran} 
\bibliography{myRef}

\begin{thebibliography}{10}
\providecommand{\url}[1]{#1}
\csname url@samestyle\endcsname
\providecommand{\newblock}{\relax}
\providecommand{\bibinfo}[2]{#2}
\providecommand{\BIBentrySTDinterwordspacing}{\spaceskip=0pt\relax}
\providecommand{\BIBentryALTinterwordstretchfactor}{4}
\providecommand{\BIBentryALTinterwordspacing}{\spaceskip=\fontdimen2\font plus
\BIBentryALTinterwordstretchfactor\fontdimen3\font minus
  \fontdimen4\font\relax}
\providecommand{\BIBforeignlanguage}[2]{{%
\expandafter\ifx\csname l@#1\endcsname\relax
\typeout{** WARNING: IEEEtran.bst: No hyphenation pattern has been}%
\typeout{** loaded for the language `#1'. Using the pattern for}%
\typeout{** the default language instead.}%
\else
\language=\csname l@#1\endcsname
\fi
#2}}
\providecommand{\BIBdecl}{\relax}
\BIBdecl

\bibitem{VehicleRoutingSurvey}
J.-F. Cordeau, G.~Laporte, M.~W. Savelsbergh, and D.~Vigo, ``Vehicle
  {R}outing,'' in \emph{Handbook in Operations Research and Management
  Science}, C.~Barnhart and G.~Laporte, Eds., no.~14.

\bibitem{TranspOnDemand}
J.-F. Cordeau, G.~Laporte, J.-Y. Potvin, and M.~W. Savelsbergh,
  ``Transportation on {D}emand,'' in \emph{Handbook in Operations Research and
  Management Science}, C.~Barnhart and G.~Laporte, Eds., no.~14.

\bibitem{VRP_book}
B.~Golden, S.~Raghavan, and E.~Wasil, Eds., \emph{The {V}ehicle {R}outing
  problem- {L}atest advances and new challenges}.\hskip 1em plus 0.5em minus
  0.4em\relax Springer, 2008.

\bibitem{rose1998deterministic}
K.~Rose, ``Deterministic annealing for clustering, compression, classification,
  regression, and related optimization problems,'' \emph{Proceedings of the
  IEEE}, vol.~86, no.~11, pp. 2210--2239, 1998.

\bibitem{xu2014aggregation}
Y.~Xu, S.~M. Salapaka, and C.~L. Beck, ``Aggregation of graph models and
  {M}arkov chains by deterministic annealing,'' \emph{Automatic Control, IEEE
  Transactions on}, vol.~59, no.~10, pp. 2807--2812, 2014.

\bibitem{kale2012maximum}
N.~V. Kale and S.~M. Salapaka, ``Maximum entropy principle-based algorithm for
  simultaneous resource location and multihop routing in multiagent networks,''
  \emph{Mobile Computing, IEEE Transactions on}, vol.~11, no.~4, pp. 591--602,
  2012.

\bibitem{salapaka2003constraints}
S.~Salapaka, A.~Khalak, and M.~Dahleh, ``Constraints on locational optimization
  problems,'' in \emph{Decision and Control, 2003. Proceedings. 42nd IEEE
  Conference on}, vol.~2.\hskip 1em plus 0.5em minus 0.4em\relax IEEE, 2003,
  pp. 1741--1746.

\bibitem{xu2014clustering}
Y.~Xu, S.~M. Salapaka, and C.~L. Beck, ``Clustering and {C}overage {C}ontrol
  for {S}ystems {W}ith {A}cceleration-{D}riven {D}ynamics,'' \emph{Automatic
  Control, IEEE Transactions on}, vol.~59, no.~5, pp. 1342--1347, 2014.

\bibitem{gray1982multiple}
R.~M. Gray and E.~D. Karnin, ``Multiple local optima in vector quantizers
  (corresp.),'' \emph{Information Theory, IEEE Transactions on}, vol.~28,
  no.~2, pp. 256--261, 1982.

\bibitem{kirkpatrick1983optimization}
S.~Kirkpatrick, C.~D. Gelatt, M.~P. Vecchi \emph{et~al.}, ``Optimization by
  simulated annealing,'' \emph{science}, vol. 220, no. 4598, pp. 671--680,
  1983.

\bibitem{braekers2011deterministic}
K.~Braekers, A.~Caris, and G.~K. Janssens, ``A deterministic annealing
  algorithm for a bi-objective full truckload vehicle routing problem in
  drayage operations,'' \emph{Procedia-Social and Behavioral Sciences},
  vol.~20, pp. 344--353, 2011.

\bibitem{braysy2008effective}
O.~Br{\"a}ysy, W.~Dullaert, G.~Hasle, D.~Mester, and M.~Gendreau, ``An
  effective multirestart deterministic annealing metaheuristic for the fleet
  size and mix vehicle-routing problem with time windows,''
  \emph{Transportation Science}, vol.~42, no.~3, pp. 371--386, 2008.

\bibitem{rose1990deterministic}
K.~Rose, ``Deterministic annealing, clustering, and optimization,'' Ph.D.
  dissertation, California Institute of Technology, 1990.

\bibitem{cover2012elements}
T.~M. Cover and J.~A. Thomas, \emph{Elements of information theory}.\hskip 1em
  plus 0.5em minus 0.4em\relax John Wiley \& Sons, 2012.

\bibitem{lloyd1982least}
S.~P. Lloyd, ``Least squares quantization in {P}{C}{M},'' \emph{Information
  Theory, IEEE Transactions on}, vol.~28, no.~2, pp. 129--137, 1982.

\bibitem{parekh2015deterministic}
P.~M. Parekh, D.~Katselis, C.~L. Beck, and S.~M. Salapaka, ``Deterministic
  annealing for clustering: {T}utorial and computational aspects,'' in
  \emph{American Control Conference (ACC), 2015}.\hskip 1em plus 0.5em minus
  0.4em\relax IEEE, 2015, pp. 2906--2911.

\bibitem{sharma2006scalable}
P.~Sharma, S.~Salapaka, and C.~Beck, ``A scalable deterministic annealing
  algorithm for resource allocation problems,'' in \emph{American Control
  Conference, 2006}.\hskip 1em plus 0.5em minus 0.4em\relax IEEE, 2006, pp.
  6--pp.

\bibitem{gillett1974heuristic}
B.~E. Gillett and L.~R. Miller, ``A heuristic algorithm for the
  vehicle-dispatch problem,'' \emph{Operations research}, vol.~22, no.~2, pp.
  340--349, 1974.

\bibitem{durbin1989analysis}
R.~Durbin, R.~Szeliski, and A.~Yuille, ``An analysis of the elastic net
  approach to the traveling salesman problem,'' \emph{Neural Computation},
  vol.~1, no.~3, pp. 348--358, 1989.

\bibitem{roehl2011maximum}
B.~A. ROEHL, ``Maximum-entropy principle approach to the multiple travelling
  salesman problem and related problems,'' Ph.D. dissertation, University of
  Illinois at Urbana-Champaign, 2011.

\bibitem{salapaka2003locational}
S.~Salapaka and A.~Khalak, ``Locational optimization problems with constraints
  on resources,'' in \emph{Proceedings of the Annual Allerton Conference on
  Communication Control and Computing}, vol.~41, no.~3.\hskip 1em plus 0.5em
  minus 0.4em\relax The University; 1998, 2003, pp. 1240--1249.

\end{thebibliography}

\end{document}